\documentclass[reqno,11pt,twoside,english]{amsart}
\usepackage[T1]{fontenc}
\usepackage[latin9]{inputenc}
\usepackage{amsthm}
\usepackage{amssymb}
\usepackage{amsmath,amsfonts,epsfig}

\makeatletter

\numberwithin{equation}{section}
\usepackage{color}
\usepackage[final,linkcolor = blue,citecolor = blue,colorlinks=true]{hyperref}

\usepackage{type1cm,ae}
\usepackage{graphicx}
\usepackage{xspace}
\usepackage{setspace}
\usepackage[english]{babel}
\usepackage{tikz}
\usetikzlibrary{arrows,patterns}

\newcommand{\rr}{{\mathbb R}}
\newcommand{\R}{{\mathbb R}}

\newcommand\wq{\infty}

\newcommand\dist{{\mathop\mathrm{\,dist\,}}}

\newcommand\De{\Delta}

\newcommand\ep{\epsilon}

\newcommand\om{{\omega}}

\newcommand\Om{{\Omega}}
\newcommand\na{{\nabla}}

\newcommand\ta{\theta}

\newcommand\pa{\partial}

\DeclareMathOperator{\divv}{div}
\numberwithin{equation}{section}
\newtheorem{theorem}{Theorem}[section]

\newtheorem{lemma}[theorem]{Lemma}

\theoremstyle{definition}
\begingroup

\endgroup

\def\XXint#1#2#3{{\setbox0=\hbox{$#1{#2#3}{\int}$}
		\vcenter{\hbox{$#2#3$}}\kern-.5\wd0}}

\makeatother

\usepackage{babel}
\begin{document}
	
	\title[Conservation law for a fourth order elliptic system in supercritical dimensions]{Conservation law for a fourth order elliptic system in supercritical dimensions}

	\author[W.-J. Qi and Z.-M. Sun]{Wen-Juan Qi and Zhao-Min Sun}
	
	\address[Wen-Juan Qi and Zhao-Min Sun]{Research Center for Mathematics and Interdisciplinary Sciences, Shandong University 266237,  Qingdao, P. R. China and Frontiers Science Center for Nonlinear Expectations, Ministry of Education, P. R. China}
	\email{wenjuan.qi@mail.sdu.edu.cn \,and\,  202217128@mail.sdu.edu.cn}
	
	%\address[Zhao-Min Sun]{Research Center for Mathematics and Interdisciplinary Sciences, Shandong University 266237,  Qingdao, P. R. China and Frontiers Science Center for Nonlinear Expectations, Ministry of Education, P. R. China}
	%\email{202217128@mail.sdu.edu.cn}

	\thanks{Corresponding author: Zhao-Min Sun.}
	
	%\thanks{Both authors are supported by the Young Scientist Program of the Ministry of Science and Technology of China (No.~2021YFA1002200) and the NSF of China (No.~12101362)}

	\begin{abstract}
		In this short note, we extend the conservation law of Lamm-Rivi\`ere [Comm. PDEs 2008] for a fourth order elliptic system to supercritical dimensions, under certain Lorentz-Sobolev integrability assumptions on the associated coefficient functions.
	\end{abstract}
	
	\maketitle
	
	{\small
		\smallskip
		\newline
		%\tableofcontents
	}
	\bigskip
	
	\arraycolsep=1pt
	
	\section{Introduction and main results}
	In his seminal work, Rivi\`ere \cite{Riviere-2007} introduced a second order linear elliptic system
	\begin{equation}\label{eq:riviere}
		-\De u=\Omega\cdot\na u \qquad \text{in } B^2,
	\end{equation}
	where $\Omega\in L^2(B^2,so(m)\otimes\Lambda^1\R^2)$. This system includes the Euler-Lagrange equation of general conformally invariant second order elliptic variational problems with quadratic growth, in particluar, harmonic mappings and prescribed mean curvature equations \cite{Riviere-2012}. Rivi\`ere successfully rewrote \eqref{eq:riviere} into divergence form, i.e. 
	$$d(\ast Adu-(\ast B)\wedge du)=0\qquad \text{in } B^2,$$ 
	provided that one found $A\in L^{\wq}\cap W^{1,2}(B^2, Gl(m))$ and $B\in W^{1,2}(B^2, M(m))$ satisfying
	$$dA-A\Omega=-d^\ast B\qquad \text{in } B^2.$$ 
	Through a fixed point argument, he succeeded in showing the exsistence of $A$ and $B$ in two dimensions, and then posed it as an open question how to find conservation laws in supercritical dimensions. Motivated by Rivi\`ere's open problem, Guo and Xiang \cite{Guo-Xiang 2 order} recently succeeded in finding a conservation law for \eqref{eq:riviere} in supercritical dimension, under certain Lorentz integrability condition on the connection matrix $\Omega$.
	
	In another pioneer work, Lamm and Rivi\`ere \cite{Lamm-Riviere-2008} extended Rivi\`ere's conservation law to a fourth order elliptic system, modeling biharmonic mappings \cite{Chang-W-Y-1999,Wang-2004-MathZ,Wang-2004-CPAM}, in the critical dimension $n=4$. More precisely, they introduced the following linear elliptic system
	\begin{equation}\label{eq:Lamm-Riviere 2008}
		\Delta^{2}u=\Delta(V\cdot\nabla u)+d^*(w\nabla u)+W\cdot\nabla u\qquad \text{in }B^4,
	\end{equation}
	where $V\in W^{1,2}(B^4,M(m)\otimes \Lambda^1\R^{4}), w\in L^{2}(B^4,M(m))$
	and $W=\nabla\omega+F$ with $\omega\in L^{2}(B^4,so(m))$ and $F\in L^{\frac{4}{3},1}(B^4,M(m)\otimes \Lambda^1\R^{4})$. 
	They proved that, if one can find $A\in W^{2,2}\cap L^{\wq}(B^4,M(m))$ and $B\in W^{1,\frac{4}{3}}(B^4,M(m)\otimes\wedge^{2}\R^{4})$ with
	\begin{equation}\label{eq:fourth order for CL}
		d\Delta A+\Delta AV-\nabla Aw+AW=d^*B\qquad \text{in }B^4,
	\end{equation}
	then system \eqref{eq:Lamm-Riviere 2008} is equivalent to the following conservation law
	\begin{equation}\label{eq:conservation law of Lamm Riviere}
		\delta[d(A \Delta u)-2 dA \Delta u+\Delta A du-A w du+dA\langle V, du\rangle-Ad\langle V,du\rangle-\ast(\ast B\wedge du)]=0.
	\end{equation}
	Lamm and Rivi\`ere \cite{Lamm-Riviere-2008}  succeeded in proving the exsistence of $A$ and $B$ on $B^4_{1/2}$, and recently, Guo et al. \cite{Guo-Xiang-Zheng-2020-fourth} refined this and sucessfully obtained the existence of $A$ and $B$ on $B^4$.
	
	Motivated by the work of Guo-Xiang \cite{Guo-Xiang 2 order}, in this paper, we would like to extend Lamm-Rivi\`ere's conservation law to supercritical dimensions. More precisely, we impose the following Lorentz-Sobolev integrability on  the coefficient functions of \eqref{eq:Lamm-Riviere 2008}: $V\in W^{1,\frac{n}{2},2}(B^n,M(m)\otimes \Lambda^1\R^{n}), w\in L^{\frac{n}{2},2}(B^n,M(m))$
	and  $W=\nabla\omega+F$
	with $\omega\in L^{\frac{n}{2},2}(B^n,so(m))$ and $F\in L^{\frac{n}{3},1}(B^n,M(m)\otimes \Lambda^1\R^{n})$. Our main result is as follows.

	% For the  fourth order elliptic system
	% \begin{equation}\label{eq:Lamm-Riviere 2008}
		% \Delta^{2}u=\Delta(V\cdot\nabla u)+d^*(w\nabla u)+W\cdot\nabla u \quad  \text{in }B^n,
		% \end{equation}
	% where $V\in W^{1,\frac{n}{2},2}(B^n,M(m)\otimes \Lambda^1\R^{n}), w\in L^{\frac{n}{2},2}(B^n,M(m))$
	% and $W$ is of the form
	% $W=\nabla\omega+F$
	% with $\omega\in L^{\frac{n}{2},2}(B^n,so(m))$ and $F\in L^{\frac{n}{3},1}(B^n,M(m)\otimes \Lambda^1\R^{n})$. 
	
	\begin{theorem}\label{thm:main theorem}
		For any $n,m\geq 4$, there exist canstants $\epsilon_{n,m},C_{n,m}>0$ satisfying the following property. Suppose $V,w,F,\omega$ satisfy the smallness condition
		\[\theta:=
		\|V\|_{W^{1,\frac{n}{2},2}}+\|w\|_{L^{\frac{n}{2},2}}+\|\omega\|_{L^{\frac{n}{2},2}}+\|F\|_{L^{\frac{n}{3},1}}<\epsilon_{n,m}.
		\]
		Then there exist $A\in W^{3,\frac{n}{3},1}\cap L^{\wq}(B^{n},Gl(m))$ and $B\in W^{1,\frac{n}{3},2}(B^{n},M(m)\otimes\wedge^{2}\R^{n})$
		such that \eqref{eq:fourth order for CL} holds in $B^n$.
		% \begin{equation}\label{eq:fourth order for CL}
			% d\Delta A+\Delta AV-\nabla Aw+AW=d^*B\quad \text{in }B^n.
			% \end{equation}
		Moreover, we have the following estimate
		\[
		\|A\|_{W^{3,\frac{n}{3},1}}+\|\dist(A,SO(m))\|_{L^{\infty}}+\|B\|_{W^{1,\frac{n}{3},2}}\leq C_{n,m}\theta.
		\]
	\end{theorem}
	
	Consequently, we obtain a conservation law: $u$ solves the fourth order Lamm-Rivi\`ere system \eqref{eq:Lamm-Riviere 2008} in $B^n$ if and only if it satisfies \eqref{eq:conservation law of Lamm Riviere} in $B^n$.
	
	\section{Auxillary results}
	In the section, we provide some necessary estimates for our later proofs.  The following result was obtained in \cite[Lemma A.3]{Lamm-Riviere-2008} when $n=4$. 
	\begin{lemma}\label{thm:4 order}
		Let $n\geq 4$ and $u\in W^{2,2}(B^n,\rr^m)$ be a solution of
		\begin{equation*}
			\begin{cases}
				\Delta^2u=*(da\wedge db)+\divv f\quad\text{in\,\,}B^n\\
				u=\frac{\partial}{\partial \nu}\Delta u=0\quad\text{on\,\,}\partial B^n\quad\text{and}\quad\int_{B^n}\Delta u=0,
			\end{cases}
		\end{equation*}
		where $a\in W^{1,\frac{n}{3},2}(B^n,\rr^m\otimes\wedge^{n-2}\rr^n)$, $b\in W^{1,n,2}(B^n,\rr^m)$ and $f\in L^{\frac{n}{3},1}(B^n,\rr^m)$. Then, 
		$$\|u\|_{L^\infty}+\|u\|_{W^{2,2}}+\|\Delta u\|_{L^{\frac{n}{2},1}}+\|\nabla\Delta u\|_{L^{\frac{n}{3},1}}\leq C\|da\|_{L^{\frac{n}{3},2}}\|db\|_{L^{n,2}}+C\|f\|_{L^{\frac{n}{3},1}}.$$
	\end{lemma}
	
	\begin{proof}
		It suffices to consider the case $n>4$. First we consider a solution $u_1$ of 
		\begin{equation*}
			\begin{cases}
				\Delta^2u_1=*(da\wedge db)\quad\text{in\,\,}B^n\\ 
				u_1=\frac{\partial}{\partial \nu}\Delta u_1=0\quad\text{on\,\,}\partial 
				B^n\quad\text{and}\quad\int_{B^n}\Delta u_1=0.
			\end{cases}
		\end{equation*}
		Since $da\wedge db\in L^{\frac{n}{4},1}\subsetneqq L^1$, the standard elliptic regularity theory implies that $u_1\in W^{4,\frac{n}{4},1}$ together with the corresponding estimate. The $L^\infty$-estimate for $u_1$, the $L^{\frac{n}{3},1}$-estimate for $\nabla\Delta u_1$ and the $L^{\frac{n}{2},1}$-estimate for $\Delta u_1$ follow from the embedding $W^{4,\frac{n}{4},1}\hookrightarrow C^0$, $W^{1,\frac{n}{4},1}\hookrightarrow L^{\frac{n}{3},1}$ and $W^{2,\frac{n}{4},1}\hookrightarrow L^{\frac{n}{2},1}$.
		
		Next we consider the solution $u_2$ of
		\begin{equation*}
			\begin{cases}
				\Delta^2u_2=\divv f\quad\text{in\,\,}B^n\\
				u_2=\frac{\partial}{\partial \nu}\Delta u_2=0\quad\text{on\,\,}\partial B^n\quad\text{and}\quad\int_{B^n}\Delta u_2=0.
			\end{cases}
		\end{equation*}
		By the $L^p$-theory for higher order elliptic system \cite[Lamma A.1]{Guo-Xiang-Zheng-2022-JMPA} and standard interpolation argument \cite{Hunt1966}, we obtain $u_2\in W^{3,\frac{n}{3},1}$. Similar to the embedding above, we can obtain the relevant estimates for $u_2$. Altogether this gives the desired estimate for $u$.
	\end{proof}
	
	The following lemma was stated in \cite[Lemma A.1]{Lamm-Riviere-2008} for $n=4$, but the proof works for all dimensions. 
	\begin{lemma}\label{thm:2 order}
		Let $n\geq 4$ and $\Tilde{f}\in L^{\frac{n}{3},2}(B^n, M(m)\otimes\wedge^1\rr^n)$. Then there exists a unique solution $v\in W^{1,\frac{n}{3},2}(B^n, M(m)\otimes\wedge^2\rr^n)$ of the Dirichlet problem
		\begin{equation*}
			\begin{cases}
				\Delta v=d\Tilde{f}\quad\text{in\,\,}B^n\\
				dv=0\quad\text{in\,\,}B^n\quad\text{and}\quad i^*_{\partial B^n}(*v)=0
			\end{cases}
		\end{equation*}
		and
		$$\|dv\|_{L^{\frac{n}{3},2}(B^n)}\leq C\|\Tilde{f}\|_{L^{\frac{n}{3},2}(B^n)}.$$
	\end{lemma}
	
%	\begin{proof}
%		This result follows from the Calder\'on-Zygmund theory and interpolation argument \cite{Hunt1966}.
%	\end{proof}
	
	\begin{lemma}\label{thm:3}
		There exists $\varepsilon>0$ such that for every $P\in W^{1,n,2}(B^n,SO(m))$  satisfying
		$$\|dP\|_{L^{n,2}(B^n)}+\|dP^{-1}\|_{L^{n,2}(B^n)}\leq \varepsilon,$$
		the only solution $C\in W^{1,\frac{n}{3},2}(B^n,M(m)\otimes\wedge^2\rr^n)$ of
		\begin{equation*}
			\begin{cases}
				d(d^*CP)=0,\\
				i^*_{\partial B^n}(*C)=0
			\end{cases}
		\end{equation*}
		is $C\equiv0$.
	\end{lemma}
	
	\begin{proof}
		Following the reduction in \cite[Proof of Lemma A.2]{Lamm-Riviere-2008}, we may assume $C=d\gamma$ for some $\gamma\in W^{1,\frac{n}{3},2}(B^n,M(m)\otimes\wedge^1\rr^n)$. 
		Since $d(d^*CP)=0$, there exists
		$D\in W^{1,\frac{n}{3},2}(B^n,M(m))$ such that $d^*CP=dD.$ This implies that
		\begin{equation*}
			\begin{cases}
				\Delta D=\pm\langle d^*C,dP\rangle\quad\text{in\,\,}B^n,\\
				D=\text{constant}\quad\qquad\text{on}\,\,\partial B^n.
			\end{cases}
		\end{equation*}
		By subtracting the constant we can assume that $D$ is zero on the boundary. Let $I_{1}$ be the first order Riesz operator with kernal $|x|^{1-n}$, and extend $C$, $P$ to $\R^n$ in a norm-bounded way. Then we can estimate
		\begin{equation}\label{eq:dD estimate}
			\begin{split}
				\|dD\|_{L^{\frac{n}{3},2}(B^n)}&\leq \|I_1*(\langle d^*C,dP\rangle)\|_{L^{\frac{n}{3},2}(\R^n)}
				\leq C\|\langle d^*C,dP\rangle\|_{L^{\frac{n}{4},1}(\R^n)}\\
				&\leq C\|dC\|_{L^{\frac{n}{3},2}(\R^n)}\|dP\|_{L^{n,2}(\R^n)}
				\leq C\varepsilon \|dC\|_{L^{n,2}(B^n)},
			\end{split}
		\end{equation}
		where we use the fact that the Riesz
		operator $I_1\colon L^{\frac{n}{4},1}(\rr^n)\to L^{\frac{n}{3},2}(\rr^n)$ is bounded.
		
		On the other hand, we have 
		\begin{equation}\label{eq:Delta C equation}
			\begin{cases}
				\Delta C=d(dDP^{-1})\quad\text{in\,\,}B^n\\
				dC=0\quad\text{in\,\,}B^n\quad\text{and}\quad i^*_{\partial B^n}(*C)=0.
			\end{cases}
		\end{equation}
		Applying Lemma~\ref{thm:2 order}, we obtain
		\begin{equation}
			\|dC\|_{L^{\frac{n}{3},2}(B^n)}\leq \tilde{C}\|dDP^{-1}\|_{L^{\frac{n}{3},2}(B^n)}\leq C\|dD\|_{L^{\frac{n}{3},2}(B^n)}.
		\end{equation}
		Combining this with \eqref{eq:dD estimate} and then choosing $\varepsilon$ sufficiently, we get that $d(*C)=0$. By \eqref{eq:Delta C equation} this implies that $*C$ is harmonic with $i^*_{\partial B^n}(*C)=0$, and thus $C\equiv0$.
	\end{proof}

	\section{Proof of Theorem~\ref{thm:main theorem}}
%	In this section, we provide the proof of Theorem~\ref{thm:main theorem}, which is very similar to \cite[Theorem 1.2]{Guo-Xiang 2 order} and \cite[Corollary 1.2]{Guo-Xiang-Zheng-2020-fourth}. 
	%Therefore, we will follow closely the statement of these papers and indicate the necessary changes when necessary.
	%The idea of the proof is almost identical to that of  \cite{Guo-Xiang-Zheng-2020-fourth} with some minor modifications.
	
	\begin{proof}[Proof of Theorem~\ref{thm:main theorem}]

		\textbf{Step 1.} Rewrite the term $W$.
		
		By standard elliptic regularity theory, there exists $\Om\in W^{1,\frac{n}{2},2}(B^{n},so(m)\otimes\wedge^{1}\R^{n})$
		such that $$d^{\ast}\Om=\om \,\,\,\text{ in }B^{n}\quad\text{ 
			with }\quad\|\Om\|_{W^{1,\frac{n}{2},2}}\le c\|\om\|_{L^{\frac{n}{2},2}}.$$
		Then for $\ep_{n,m}$ is sufficiently small with $\|\om\|_{L^{\frac{n}{2},2}}\le\ep_{n,m}$, we can find $U\in W^{2,\frac{n}{2},2}(B^{n},so(m))$, $P=e^{U}\in W^{2,\frac{n}{2},2}(B^{n},SO(m))$ and $\xi\in W^{2,\frac{n}{2},2}(B^{n},so(m)\otimes\wedge^{2}\R^{n})$ such that
		$$\Om=P^{-1}dP+P^{-1}d^{\ast}\xi P \quad\text{in}\,\,\, B^{n} \qquad\text{and}\qquad d\left(i_{\pa B^{n}}^{\ast}\ast\xi\right)=0 \quad\text{on}\,\,\, \pa B^{n}$$
		with
		$$\|P\|_{W^{2,\frac{n}{2},2}}+\|\xi\|_{W^{2,\frac{n}{2},2}}\le C_{n,m}\|\Om\|_{W^{1,\frac{n}{2},2}}\le C_{n,m}\ep_{n,m}.$$
		The proof for this is completely similar to that of \cite[Theorem A.5]{Lamm-Riviere-2008}, thus we omit it here.
		
		Direct calculation yields
		\begin{equation*}
			W=d\omega+F=d(d^{\ast}\Omega)+F=-P^{-1}d\Delta P+K_{1},\label{eq: 2.5}
		\end{equation*}
		where 
		$
		K_{1}=-dP^{-1}\Delta P-d\langle dP^{-1},dP\rangle-d\langle P^{-1}d^{\ast}\xi,dP\rangle-d\langle dP^{-1},d^{\ast}\xi P\rangle+F.
		$
		
		It follows from the Sobolev embedding $W^{1,\frac{n}{2},2}\subset L^{n,2}
		$ that $K_{1}\in L^{\frac{n}{3},1}(B^{n})$ with 
		\[
		\|K_{1}\|_{L^{\frac{n}{3},1}(B^{n})}\leq C\|\omega\|_{L^{\frac{n}{2},2}(B^{n})}+C\|F\|_{L^{\frac{n}{3},1}(B^{n})}\leq C_{n,m}\ta.
		\]

		% \textbf{Step 3.} Reduce to an equivalent problem 
		\textbf{Step 2.} Reduce to an equivalent problem. 
		
		Note that there exists a pair $(A,B)$ solves \eqref{eq:fourth order for CL} if and only if for $(\Tilde{A},B)$ with $\tilde{A}=AP^{-1}$, we have
		\begin{equation}\label{eq: 2.8}
			d\De\tilde{A}+\De\tilde{A}K_{2}+\na^{2}\tilde{A}K_{3}+d\tilde{A}K_{4}+\tilde{A}K_{5}=d^{\ast}BP^{-1}\quad \text{ in }\,\,B^n,
		\end{equation}
		where 
		\begin{equation*}
			\begin{split}
				&K_2=PV P^{-1}+\na PP^{-1},\,\,K_3=2\na P P^{-1},\\
				&K_4=2\na P V P^{-1}-Pw P^{-1}+2\na^2 PP^{-1}+\De PP^{-1}, \\
				&K_5=\De PV P^{-1}-\na Pw P^{-1}+PK_1P^{-1}
			\end{split}
		\end{equation*}
		with the estimate
		\begin{equation}
			\|K_{2}\|_{W^{1,\frac{n}{2},2}(B^{n})}+\|K_{3}\|_{W^{1,\frac{n}{2},2}(B^{n})}+\|K_{4}\|_{L^{\frac{n}{2},2}(B^{n})}+\|K_{5}\|_{L^{\frac{n}{3},1}(B^{n})}<C\ta.\label{eq: 2.9}
		\end{equation}
		
		We now use an extension argument from \cite{Guo-Xiang 2 order, Guo-Xiang-Zheng-2020-fourth} as follows: extend $K_{i}$ ($2\le i\le5$), $U$ and $P$ to the whole space $\mathbb{R}^n$ with compact support in $B_{2}^{n}$ in a norm-bounded way. For simplicity, we will still use the original notations for representation. Our strategy is to solve \eqref{eq: 2.8} in
		the enlarged region $B_2^n$.

		\textbf{Claim.}\label{claim:AB}
		Set $\tilde{A}=\bar{A}+I$, there exsits $(\Bar{A},B)$ such that
		\begin{equation}
			\begin{cases}
				\De^{2}\bar{A}=d^{\ast}\left(\De\bar{A}{K}_{2}+\na^{2}\bar{A}{K}_{3}+d\bar{A}{K}_{4}+\bar{A}{K}_{5}+{K}_{5}-d^{\ast}BP^{-1}\right)  \quad\text{in }B_{2}^{n},\\
				\De B=d\left[\left(d\De\bar{A}+\De\bar{A}{K}_{2}+\na^{2}\bar{A}{K}_{3}+d\bar{A}{K}_{4}+\bar{A}{K}_{5}+{K}_{5}\right)P\right] \quad\,\, \text{in }B_{2}^{n},\\
				\bar{A}=\frac{\pa\De\bar{A}}{\pa\nu}=0  \quad\text{on }\,\,\pa B_{2}^{n}\qquad \text{and} \qquad \int_{B_{2}^{n}}\De\bar{A}=0, \\
				dB=0  \quad\text{ in }\,\,B_{2}^{n} \qquad \text{and}\qquad i_{\pa B_{2}^{n}}^{\ast}(\ast B)=0.
			\end{cases}\label{eq: 2.11-1}
		\end{equation}
		Moreover, we have
		\begin{equation}\label{eq:claim estimate}
			\|\Bar{A}\|_{W^{3,\frac{n}{3},1}(B^{n})}+\|\Bar{A}\|_{L^{\wq}(B^{n})}+\|B\|_{W^{1,\frac{n}{3},2}(B^{n})}\leq C\theta\leq \epsilon_{n,m}.
		\end{equation}

		\textbf{Proof of Claim.}
		% Also, we will use the fixed point argument as in \cite{Lamm-Riviere-2008}. First, we introduce 
		Similiar to \cite{Guo-Xiang 2 order, Guo-Xiang-Zheng-2020-fourth}, we introduce a Banach space $\mathbb{H}=(\mathbb{H},\|\cdot\|_{\mathbb{H}})$:
		\[
		\mathbb{H}=\left\{ (u,v)\in W^{2,\frac{n}{2},2}\cap L^{\wq}(B_{2}^{n},M(m))\times W^{1,\frac{n}{3},2}(B_{2}^{n},M(m)\otimes\wedge^{2}\R^{n}):\|(u,v)\|_\mathbb{H}\leq 1
		\right\} 
		% satisfies }\eqref{eq: 2.11-2}\right\} 
	\]
	where
	$
	\|(u,v)\|_{\mathbb{H}}\equiv\|u\|_{W^{2,\frac{n}{2},2}(B^{n})}+\|u\|_{L^{\wq}(B^{n})}+\|v\|_{W^{1,\frac{n}{3},2}(B^{n})}.
	$
	
	Standard elliptic regularity theory implies that for each $(u,v)\in\mathbb{H}$, there exists a
	unique solution $\bar{u}\in W^{2,\frac{n}{2},2}(B^{n}_2)$
	satisfying 
	\begin{equation}\label{eq:ubar}
		\begin{cases}
			\De^{2}\bar{u}=d^{\ast}\left(\De u{K}_{2}+\na^{2}u{K}_{3}+duK_{4}+uK_{5}-d^{\ast}vP^{-1}+{K}_{5}\right) \quad \text{in\,\,}B_{2}^{n},\\
			\bar{u}=\frac{\pa\De\bar{u}}{\pa\nu}=0  \quad \text{on}\,\,\,\,\pa B_{2}^{n}\qquad \text{and}\qquad\int_{B_{2}^{n}}\De\bar{u}=0.
		\end{cases}
	\end{equation}
	% \colon=d^{\ast}\left(f-d^{\ast}v\tilde{P}^{-1}\right)
	Set
	$$f:=\De u{K}_{2}+\na^{2}u{K}_{3}+du{K}_{4}+u{K}_{5}+{K}_{5}.$$
	It follows from $u\in W^{2,\frac{n}{2},2}$ and $v\in W^{1,\frac{n}{3},2}$ that $f\in L^{\frac{n}{3},1}(B_{2}^{n})$ with
	$$
	\|f\|_{L^{\frac{n}{3},1}}\le c\ta\left(\|u\|_{W^{2,\frac{n}{2},2}}+\|u\|_{L^{\wq}}+1\right).
	$$
	%Note that $u\in W^{2,\frac{n}{2},2}$ and $v\in W^{1,\frac{n}{3},2}$ implies that $f\in L^{\frac{n}{3},1}(B_{2}^{n})$, where
	%$$f:=\De u{K}_{2}+\na^{2}u{K}_{3}+du{K}_{4}+u{K}_{5}+{K}_{5}$$
	%with $$\|f\|_{L^{\frac{n}{3},1}}\le c\ta\left(\|u\|_{W^{2,\frac{n}{2},2}}+\|u\|_{L^{\wq}}+1\right).$$
	Together with
	\[
	d^{\ast}(d^{\ast}v\tilde{P}^{-1})=-\langle d^{\ast}v,d\tilde{P}^{-1}\rangle=\pm\ast(d\ast v\wedge d\tilde{P}^{-1}),
	\]
	we can apply Lemma \ref{thm:4 order} to deduce that $\bar{u}\in W^{3,\frac{n}{3},1}(B_{2}^{n})$
	with 
	\[
	\|\bar{u}\|_{L^{\wq}}+\|\bar{u}\|_{W^{2,\frac{n}{2},2}}+\|\De\bar{u}\|_{L^{\frac{n}{2},1}}+\|d\De\bar{u}\|_{L^{\frac{n}{3},1}}\le C\ta(\|u\|_{W^{2,\frac{n}{2},2}}+\|u\|_{L^{\wq}}+1+\|v\|_{W^{1,\frac{n}{3},2}}).
	\]
	Lemma \ref{thm:2 order} implies that there exists a unique $\bar{v}\in W^{1,\frac{n}{3},2}(B_{2}^{n},M(m)\otimes\wedge^{2}\R^{n})$ solving 
	\begin{equation}\label{eq:vbar}
		\begin{cases}
			\De\bar{v}=d\left[\left(d\De\bar{u}+\De\bar{u}{K}_{2}+\na^{2}\bar{u}{K}_{3}+d\bar{u}{K}_{4}+\bar{u}{K}_{5}+{K}_{5}\right)P\right]  &\text{in }B_{2}^{n},\\
			d\bar{v}=0  \quad\text{in }\,\, B_{2}^{n}\qquad\text{and}\qquad
			i_{\pa B_{2}^{n}}^{\ast}(\ast\bar{v})=0.
		\end{cases}
	\end{equation}
	where $\tilde{f}=\left(d\De\bar{u}+\De\bar{u}{K}_{2}+\na^{2}\bar{u}{K}_{3}+d\bar{u}{K}_{4}+\bar{u}{K}_{5}+{K}_{5}\right)P.$ Furthermore, we have
	\[
	\|d\bar{v}\|_{L^{\frac{n}{3},2}}\le\|\Tilde{f}\|_{L^{\frac{n}{3},2}}\le \|\tilde{f}\|_{L^{\frac{n}{3},1}}\le C\ta\left(\|u\|_{W^{2,\frac{n}{2},2}}+\|u\|_{L^{\wq}}+\|v\|_{W^{1,\frac{n}{3},2}}+1\right).
	\]
	Remind that all the norms in the above are taken over $B_{2}^{n}$, thus we omit them for brevity.
	Hence, for any $(u,v)\in\mathbb{H}$, there exists a unique $(\bar{u},\bar{v})\in\mathbb{H}$
	solving systems \eqref{eq:ubar} and \eqref{eq:vbar}.
	% \[
	% \begin{cases}
		% \De^{2}\bar{u}=d^{\ast}\left(\De u\tilde{K}_{2}+\na^{2}u\tilde{K}_{3}+du\tilde{K}_{4}+u\tilde{K}_{5}-d^{\ast}v\tilde{P}^{-1}+\tilde{K}_{5}\right) & \text{in }B_{2}^{n},\\
		% \De\bar{v}=d\left[\left(d\De\bar{u}+\De\bar{u}\tilde{K}_{2}+\na^{2}\bar{u}\tilde{K}_{3}+d\bar{u}\tilde{K}_{4}+\bar{u}\tilde{K}_{5}+\tilde{K}_{5}\right)\tilde{P}\right] & \text{in }B_{2}^{n},\\
		% \bar{u}=\frac{\pa\De\bar{u}}{\pa\nu}=0  \quad \text{on}\,\,\,\,\pa B_{2}^{n}\qquad \text{and}\qquad\int_{B_{2}^{n}}\De\bar{u}=0,\\
		% d\bar{v}=0  \quad\text{in }\,\, B_{2}^{n}\qquad\text{and}\qquad
		% i_{\pa B_{2}^{n}}^{\ast}(\ast\bar{v})=0.
		% \end{cases}
	% \]
	Moreover, there
	exists $C>0$ such that
	\begin{equation}
		\|(\bar{u},\bar{v})\|_{\mathbb{H}}\le C\ta(\|(u,v)\|_{\mathbb{H}}+1)\label{eq: apriori estimate}
	\end{equation}
	This implies if $\ep_{n,m}$ is sufficiently small such that $C\ta\leq\frac{1}{2}$, then $(\bar{u},\bar{v})\in\mathbb{H}.$ If we define the mapping $T\colon \mathbb{H}\to\mathbb{H}$ by
	$$T(u,v)=(\bar{u},\bar{v}),$$
	then $T$ is a contraction operator on $\mathbb{H}$. Indeed, we can use a similar argument to obtain 
	\[
	\|T(u_{1},v_{1})-T(u_{2},v_{2})\|_{\mathbb{H}}\le\frac{1}{2}\|(u_{1},v_{1})-(u_{2},v_{2})\|_{\mathbb{H}}.
	\]
	By the standard fixed point theorem, there exsit a unique $(\bar{A},B)\in\mathbb{H}$ such that
	\[
	T(\bar{A},B)=(\bar{A},B).
	\]
	That is, $(\bar{A},B)$ solves problem \eqref{eq: 2.11-1}, where $\bar{A}\in W^{3,\frac{n}{3},1}(B_{2}^{n})$ and $B\in W^{1,\frac{n}{3},2}(B_{2}^{n})$. Moreover,
	by \eqref{eq: apriori estimate}, we obtain the estimate \eqref{eq:claim estimate}. The proof of Claim is complete.
	%$$\|(\bar{A},B)\|_{\mathbb{H}}\leq C\ta\leq\ep_{n,m}.$$

	\textbf{Step 3.} Solve the orginal system.
	% \
	% \newline
	% \noindent
	
	By Claim, 
	% now we set
	% $\tilde{A}=\bar{A}+I$, 
	$\tilde{A}=\bar{A}+I\in W^{3,\frac{n}{3},1}(B_{2}^{n})$, $B\in W^{1,\frac{n}{3},2}(B_{2}^{n})$
	and they satisfy 
	\[
	d^{\ast}\left(d\De\tilde{A}+\De\tilde{A}K_{2}+\na^{2}\tilde{A}K_{3}+d\tilde{A}K_{4}+\tilde{A}K_{5}-d^{\ast}BP^{-1}\right)=0\qquad\text{in }B_{2}^{n}.
	\]
	It follows from the Hodge decomposition that 
	% [Theorem 2.4.14]
	\[
	d\De\tilde{A}+\De\tilde{A}K_{2}+\na^{2}\tilde{A}K_{3}+d\tilde{A}K_{4}+\tilde{A}K_{5}-d^{\ast}BP^{-1}=d^{\ast}C+h,
	\]
	where $C\in W^{1,\frac{n}{3},2}(B_{2}^{n},M(m)\otimes\wedge^{2}\R^{n})$
	satisfies $i_{\pa B_{2}^{n}}^{\ast}(\ast C)=0$ and $h$ is a harmonic
	$1$-form in $B_{2}^{n}$. 
	By the same argument as that of \cite[Page 10]{Guo-Xiang-Zheng-2020-fourth}, we deduce that $h\equiv0$. Hence, we have 
	\[
	d\De\tilde{A}+\De\tilde{A}K_{2}+\na^{2}\tilde{A}K_{3}+d\tilde{A}K_{4}+\tilde{A}K_{5}-d^{\ast}BP^{-1}=d^{\ast}C
	\]
	with $i_{\pa B_{2}^{n}}^{\ast}(\ast C)=0$. By directly calculation, we have $d(d^*CP)=0$. Applying Lemma~\ref{thm:3}, we infer that $C\equiv0$ by choosing $\ep_{m}$ sufficiently small. The proof is thus complete. 
\end{proof}

%\subsection*{A. Appendix}
%\appendix
%\renewcommand{\thelemma}{A.\arabic{lemma}}
%\renewcommand{\theequation}{A.\arabic{equation}}

\subsection*{Acknowledgements}\small
Both authors are supported by the Young Scientist Program of the Ministry of Science and Technology of China (No.~2021YFA1002200), the NSF of China (No.~12101362) and the Taishan Scholar Program and the Natural Science Foundation of Shandong Province (No. ZR2022YQ01, ZR2021QA003). The authors would like to thank their supervisor Prof.~Chang-Yu Guo for posing this question to them and for many useful conservations.

\end{document}